\font\fifteenrm=cmr10 scaled\magstep2 
\font\fifteeni=cmmi10 scaled\magstep2
\font\fifteensy=cmsy10 scaled\magstep2
\font\fifteenbf=cmbx10 scaled\magstep2
\font\fifteentt=cmtt10 scaled\magstep2
\font\fifteenit=cmti10 scaled\magstep2
\font\fifteensl=cmsl10 scaled\magstep2
\font\fifteenam=msam10 scaled\magstep2
\font\fifteenbm=msbm10 scaled\magstep2
\font\fifteenex=cmex10 scaled\magstep2
\font\fifteensc=cmcsc10 scaled\magstep2 
\font\twelverm=cmr10 at 12pt
\font\twelvei=cmmi10 at 12pt
\font\twelvesy=cmsy10 at 12pt
\font\twelvebf=cmbx10 at 12pt
\font\twelvett=cmtt10 at 12pt
\font\twelveit=cmti10 at 12pt
\font\twelvesl=cmsl10 at 12pt
\font\twelveam=msam10 at 12pt
\font\twelvebm=msbm10 at 12pt
\font\twelveex=cmex10 at 12pt
\font\twelvesc=cmcsc10 at 12pt
\font\elevenrm=cmr10 scaled\magstephalf 
\font\eleveni=cmmi10 scaled\magstephalf
\font\elevensy=cmsy10 scaled\magstephalf
\font\elevenbf=cmbx10 scaled\magstephalf
\font\eleventt=cmtt10 scaled\magstephalf
\font\elevenit=cmti10 scaled\magstephalf
\font\elevensl=cmsl10 scaled\magstephalf
\font\elevenam=msam10 scaled\magstephalf
\font\elevenbm=msbm10 scaled\magstephalf
\font\elevenex=cmex10 scaled\magstephalf
\font\elevensc=cmcsc10 scaled\magstephalf
\font\tenrm=cmr10
\font\teni=cmmi10
\font\tensy=cmsy10
\font\tenbf=cmbx10
\font\tentt=cmtt10
\font\tenit=cmti10
\font\tensl=cmsl10
\font\tenam=msam10
\font\tenbm=msbm10
\font\tenex=cmex10
\font\tensc=cmcsc10
\font\ninerm=cmr9
\font\ninei=cmmi9
\font\ninesy=cmsy9
\font\ninebf=cmbx9
\font\ninett=cmtt9
\font\nineit=cmti9
\font\ninesl=cmsl9
\font\nineam=msam9
\font\ninebm=msbm9
\font\nineex=cmex9
\font\ninesc=cmcsc9
\font\eightrm=cmr8
\font\eighti=cmmi8
\font\eightsy=cmsy8
\font\eightbf=cmbx8
\font\eighttt=cmtt8
\font\eightit=cmti8
\font\eightsl=cmsl8
\font\eightam=msam8
\font\eightbm=msbm8
\font\eightex=cmex8
\font\eightsc=cmcsc8
\font\sevenrm=cmr7
\font\seveni=cmmi7
\font\sevensy=cmsy7
\font\sevenbf=cmbx7

\font\sevenam=msam7
\font\sevenbm=msbm7

\font\sixrm=cmr6
\font\sixi=cmmi6
\font\sixsy=cmsy6

\font\sixam=msam6
\font\sixbm=msbm6

\font\fiverm=cmr5
\font\fivei=cmmi5
\font\fivesy=cmsy5

\font\fiveam=msam5
\font\fivebm=msbm5

\font\fourrm=cmr5 at 4pt
\font\fouri=cmmi5 at 4pt
\font\foursy=cmsy5 at 4pt

\font\fouram=msam5 at 4pt
\font\fourbm=msbm5 at 4pt

\skewchar\twelvei='177 \skewchar\eleveni='177\skewchar\teni='177
\skewchar\ninei='177 \skewchar\eighti='177\skewchar\seveni='177 
\skewchar\sixi='177 \skewchar\fivei='177 \skewchar\fouri='177
\skewchar\twelvesy='60 \skewchar\elevensy='60 \skewchar\tensy='60
\skewchar\ninesy='60 \skewchar\eightsy='60 \skewchar\sevensy='60 
\skewchar\sixsy='60 \skewchar\fivesy='60 \skewchar\foursy='60
\newfam\itfam
\newfam\slfam
\newfam\bffam
\newfam\ttfam
\newfam\scfam
\newfam\amfam
\newfam\bmfam
\def\eightbig#1{{\hbox{$\left#1\vbox to 6.5pt{}\voidright $}}}
\def\eightBig#1{{\hbox{$\left#1\vbox to 7.5pt{}\voidright $}}}
\def\eightbigg#1{{\hbox{$\left#1\vbox to 10pt{}\voidright $}}}
\def\eightBigg#1{{\hbox{$\left#1\vbox to 13pt{}\voidright $}}}
\def\ninebig#1{{\hbox{$\left#1\vbox to 7.5pt{}\voidright $}}}
\def\nineBig#1{{\hbox{$\left#1\vbox to 8.5pt{}\voidright $}}}
\def\ninebigg#1{{\hbox{$\left#1\vbox to 11.5pt{}\voidright $}}}
\def\nineBigg#1{{\hbox{$\left#1\vbox to 14.5pt{}\voidright $}}}
\def\tenbig#1{{\hbox{$\left#1\vbox to 8.5pt{}\voidright $}}}
\def\tenBig#1{{\hbox{$\left#1\vbox to 9.5pt{}\voidright $}}}
\def\tenbigg#1{{\hbox{$\left#1\vbox to 12.5pt{}\voidright $}}}
\def\tenBigg#1{{\hbox{$\left#1\vbox to 16pt{}\voidright $}}}
\def\elevenbig#1{{\hbox{$\left#1\vbox to 9pt{}\voidright $}}}
\def\elevenBig#1{{\hbox{$\left#1\vbox to 10.5pt{}\voidright $}}}
\def\elevenbigg#1{{\hbox{$\left#1\vbox to 14pt{}\voidright $}}}
\def\elevenBigg#1{{\hbox{$\left#1\vbox to 17.5pt{}\voidright $}}}
\def\twelvebig#1{{\hbox{$\left#1\vbox to 10pt{}\voidright $}}}
\def\twelveBig#1{{\hbox{$\left#1\vbox to 11pt{}\voidright $}}}
\def\twelvebigg#1{{\hbox{$\left#1\vbox to 15pt{}\voidright $}}}
\def\twelveBigg#1{{\hbox{$\left#1\vbox to 19pt{}\voidright $}}}
\def\fifteenbig#1{{\hbox{$\left#1\vbox to 12pt{}\voidright $}}}
\def\fifteenBig#1{{\hbox{$\left#1\vbox to 13.5pt{}\voidright $}}}
\def\fifteenbigg#1{{\hbox{$\left#1\vbox to 18pt{}\voidright $}}}
\def\fifteenBigg#1{{\hbox{$\left#1\vbox to 23pt{}\voidright $}}}
\def\voidright{\right.\nulldelimiterspace=0pt \mathsurround=0pt }
\def\fifteenpoint{
  \textfont0=\fifteenrm \scriptfont0=\twelverm \scriptscriptfont0=\tenrm
  \def\rm{\fam0 \fifteenrm}%
  \textfont1=\fifteeni \scriptfont1=\twelvei \scriptscriptfont1=\teni
  \textfont2=\fifteensy \scriptfont2=\twelvesy \scriptscriptfont2=\tensy
  \textfont3=\fifteenex \scriptfont3=\fifteenex \scriptscriptfont3=\fifteenex
  \def\it{\fam\itfam\fifteenit}\textfont\itfam=\fifteenit
  \def\sl{\fam\slfam\fifteensl}\textfont\slfam=\fifteensl
  \def\bf{\fam\bffam\fifteenbf}\textfont\bffam=\fifteenbf 
    \scriptfont\bffam=\twelvebf\scriptscriptfont\bffam=\tenbf
  \def\tt{\fam\ttfam\fifteentt}\textfont\ttfam=\fifteentt
  \def\sc{\fam\scfam\fifteensc}\textfont\scfam=\fifteensc
  \def\am{\fam\amfam\fifteenam}\textfont\amfam=\fifteenam
    \scriptfont\amfam=\twelveam\scriptscriptfont\amfam=\tenam
  \def\bm{\fam\bmfam\fifteenbm}\textfont\bmfam=\fifteenbm
    \scriptfont\bmfam=\twelvebm\scriptscriptfont\bmfam=\tenbm
  \baselineskip=21pt \rm
  \let\big=\fifteenbig\let\Big=\fifteenBig\let\bigg=\fifteenbigg
  \let\Bigg=\fifteenBigg}
\def\twelvepoint{
  \textfont0=\twelverm \scriptfont0=\ninerm \scriptscriptfont0=\sevenrm
  \def\rm{\fam0 \twelverm}%
  \textfont1=\twelvei \scriptfont1=\ninei \scriptscriptfont1=\seveni
  \textfont2=\twelvesy \scriptfont2=\ninesy \scriptscriptfont2=\sevensy
  \textfont3=\twelveex \scriptfont3=\twelveex \scriptscriptfont3=\twelveex
  \def\it{\fam\itfam\twelveit}\textfont\itfam=\twelveit
  \def\sl{\fam\slfam\twelvesl}\textfont\slfam=\twelvesl
  \def\bf{\fam\bffam\twelvebf}\textfont\bffam=\twelvebf 
    \scriptfont\bffam=\ninebf\scriptscriptfont\bffam=\sevenbf
  \def\tt{\fam\ttfam\twelvett}\textfont\ttfam=\twelvett
  \def\sc{\fam\scfam\twelvesc}\textfont\scfam=\twelvesc
  \def\am{\fam\amfam\twelveam}\textfont\amfam=\twelveam
    \scriptfont\amfam=\nineam\scriptscriptfont\amfam=\sevenam
  \def\bm{\fam\bmfam\twelvebm}\textfont\bmfam=\twelvebm
    \scriptfont\bmfam=\ninebm\scriptscriptfont\bmfam=\sevenbm
  \baselineskip=17.8pt \rm 
  \def\looselineskip{\baselineskip=18.5pt plus 1.8pt}%
  \def\tightlineskip{\baselineskip=16.5pt}%
  \def\verytightlineskip{\baselineskip=15pt}%
  \let\big=\twelvebig\let\Big=\twelveBig\let\bigg=\twelvebigg
  \let\Bigg=\twelveBigg  }
\def\elevenpoint{
  \textfont0=\elevenrm \scriptfont0=\ninerm \scriptscriptfont0=\sixrm
  \def\rm{\fam0 \elevenrm}%
  \textfont1=\eleveni \scriptfont1=\ninei \scriptscriptfont1=\sixi
  \textfont2=\elevensy \scriptfont2=\ninesy \scriptfont2=\sixsy 
  \textfont3=\elevenex \scriptfont3=\elevenex \scriptfont3=\elevenex
  \def\it{\fam\itfam\elevenit}\textfont\itfam=\elevenit
  \def\sl{\fam\slfam\elevensl}\textfont\slfam=\elevensl
  \def\bf{\fam\bffam\elevenbf}\textfont\bffam=\elevenbf
  \def\tt{\fam\ttfam\eleventt}\textfont\ttfam=\eleventt
  \def\sc{\fam\scfam\elevensc}\textfont\scfam=\elevensc
  \def\am{\fam\amfam\elevenam}\textfont\amfam=\elevenam
    \scriptfont\amfam=\nineam\scriptscriptfont\amfam=\sixam
  \def\bm{\fam\bmfam\elevenbm}\textfont\bmfam=\elevenbm
    \scriptfont\bmfam=\ninebm\scriptscriptfont\bmfam=\sixbm
  \baselineskip=15.1pt \rm
  \def\looselineskip{\baselineskip=16pt plus 1.5pt}%
  \def\tightlineskip{\baselineskip=14pt}%
  \def\verytightlineskip{\baselineskip=13pt}%
  \let\big=\elevenbig\let\Big=\elevenBig\let\bigg=\elevenbigg
  \let\Bigg=\elevenBigg  }
\def\tenpoint{
  \textfont0=\tenrm \scriptfont0=\eightrm \scriptscriptfont0=\fiverm
  \def\rm{\fam0 \tenrm}%
  \textfont1=\teni \scriptfont1=\eighti \scriptscriptfont1=\fivei
  \textfont2=\tensy \scriptfont2=\eightsy \scriptfont2=\fivesy 
  \textfont3=\tenex \scriptfont3=\tenex \scriptfont3=\tenex
  \def\it{\fam\itfam\tenit}\textfont\itfam=\tenit
  \def\sl{\fam\slfam\tensl}\textfont\slfam=\tensl
  \def\bf{\fam\bffam\tenbf}\textfont\bffam=\tenbf
  \def\tt{\fam\ttfam\tentt}\textfont\ttfam=\tentt
  \def\sc{\fam\scfam\tensc}\textfont\scfam=\tensc
  \def\am{\fam\amfam\tenam}\textfont\amfam=\tenam
    \scriptfont\amfam=\eightam \scriptscriptfont\amfam=\fiveam
  \def\bm{\fam\bmfam\tenbm}\textfont\bmfam=\tenbm
    \scriptfont\bmfam=\eightbm \scriptscriptfont\bmfam=\fivebm
  \baselineskip=14pt \rm
  \def\looselineskip{\baselineskip=14.8pt plus1.5pt}
  \def\tightlineskip{\baselineskip=13.6pt}%
  \def\verytightlineskip{\baselineskip=13pt}%
  \let\big=\tenbig\let\Big=\tenBig\let\bigg=\tenbigg\let\Bigg=\tenBigg  }
\def\ninepoint{
  \textfont0=\ninerm \scriptfont0=\sevenrm \scriptscriptfont0=\fourrm
  \def\rm{\fam0 \ninerm}%
  \textfont1=\ninei \scriptfont1=\seveni \scriptscriptfont1=\fouri
  \textfont2=\ninesy \scriptfont2=\sevensy \scriptfont2=\foursy 
  \textfont3=\nineex \scriptfont3=\nineex \scriptfont3=\nineex
  \def\it{\fam\itfam\nineit}\textfont\itfam=\nineit
  \def\sl{\fam\slfam\ninesl}\textfont\slfam=\ninesl
  \def\bf{\fam\bffam\ninebf}\textfont\bffam=\ninebf
  \def\tt{\fam\ttfam\ninett}\textfont\ttfam=\ninett
  \def\sc{\fam\scfam\ninesc}\textfont\scfam=\ninesc
  \def\am{\fam\amfam\nineam}\textfont\amfam=\nineam
    \scriptfont\amfam=\nineam\scriptscriptfont\amfam=\fouram
  \def\bm{\fam\bmfam\ninebm}\textfont\bmfam=\ninebm
    \scriptfont\bmfam=\ninebm\scriptscriptfont\bmfam=\fourbm
  \baselineskip=12.6pt \rm
  \let\big=\ninebig\let\Big=\nineBig\let\bigg=\ninebigg
  \let\Bigg=\nineBigg  }
\def\eightpoint{
  \textfont0=\eightrm \scriptfont0=\fiverm \scriptscriptfont0=\fourrm
  \def\rm{\fam0 \eightrm}%
  \textfont1=\eighti \scriptfont1=\fivei \scriptscriptfont1=\fouri
  \textfont2=\eightsy \scriptfont2=\fivesy \scriptfont2=\foursy 
  \textfont3=\eightex \scriptfont3=\eightex \scriptfont3=\eightex
  \def\it{\fam\itfam\eightit}\textfont\itfam=\eightit
  \def\sl{\fam\slfam\eightsl}\textfont\slfam=\eightsl
  \def\bf{\fam\bffam\eightbf}\textfont\bffam=\eightbf
  \def\tt{\fam\ttfam\eighttt}\textfont\ttfam=\eighttt
  \def\sc{\fam\scfam\eightsc}\textfont\scfam=\eightsc
  \def\am{\fam\amfam\eightam}\textfont\amfam=\eightam
    \scriptfont\amfam=\eightam\scriptscriptfont\amfam=\fouram
  \def\bm{\fam\bmfam\eightbm}\textfont\bmfam=\eightbm
    \scriptfont\bmfam=\eightbm\scriptscriptfont\bmfam=\fourbm
  \baselineskip=11.2pt \rm
  \let\big=\eightbig\let\Big=\eightBig\let\bigg=\eightbigg
  \let\Bigg=\eightBigg  }

\twelvepoint
\nopagenumbers
\hsize=6in\vsize=8.8in

\parskip=1pt plus 1pt

\newif\ifSpecialhead\Specialheadfalse
\newbox\specialheadbox

\def\specialhead #1\par{\Specialheadtrue\setbox\specialheadbox=\hbox{#1}}
\headline={{\ifSpecialhead\box\specialheadbox\global\Specialheadfalse\else
     \ifnum\pageno<0{\hfill\quad{\twelvebf\folio}}%
     \else\ifnum\pageno<2\hfill
     \else\hfill\twelvepoint\sc\firstmark\quad{\twelvebf\folio}\fi\fi\fi}}

\def\title#1\par{\bigskip{\def\cr{\par\center}\center\fifteenbf #1\par}\medskip}
\def\subtitle#1\par{\centerline{\fifteenrm #1}\medskip}
\def\author#1\par{\medskip{\def\cr{\par\center\twelvesc}\fifteensc\center#1\par}}
\def\center#1\par{\hfil #1\hfil\par}
\def\abstract.#1\par{\message{Abstract.}%
                    \medskip{\narrower\narrower\tenpoint\tightlineskip
                        \noindent{\bf Abstract.}#1\par}\medskip\noindent}
\def\bigabstract.#1\par{\message{Abstract.}%
                         \medskip{\narrower\narrower\tightlineskip
                         \noindent{\bf Abstract. }#1\par}\medskip\noindent}
\def\acknowledgement#1\par{\footnote{}{#1}}
\def\sectionskip{\Goodbreak\vskip 25pt plus 15pt minus 5pt}
\def\secnumber{\ifquiet
               \else\ifNoSections
                    \else\sectionsymbol\the\secno\quad\fi\fi}
\def\section#1\par{ \NoSectionsfalse\par\sectionskip\proofdepth=0\claimno=0
 \ifquiet\else\advance\secno by1\fi\toks0={#1}
 \immediate\write16{\ifquiet\else Section \the\secno\space\fi
                    \the\toks0}%
 \mark{\secnumber #1}%
 {\fifteenpoint\bf\noindent\secnumber #1}\nobreak\bigskip\quietoff
 \nobreak\noindent}
\def\quiet{\quiettrue}
\def\QUIET{\QUIETtrue\quiettrue}

\def\quietoff{\ifQUIET\else\quietfalse\fi}
\newif\ifquiet
\newif\ifQUIET
\newif\ifNoSections
\newcount\claimtype
\newcount\secno
\newcount\claimno
\newcount\subclaimno
\newcount\subsubclaimno
\newcount\subsubsubclaimno
\newcount\proofdepth
\def\subclaimnumber{\ifquiet\else\ifcase\subclaimno\or A\or B\or C\or D\or E\or
     F\or G\or H\or I\or J\or K\or L\or M\or N\or O\or P\fi\fi}
\def\subsubclaimnumber{\ifquiet\else\ifcase\subsubclaimno\or i\or ii\or iii\or 
   iv\or v\or vi\or vii\or viii\or ix\or x\or xi\or xii\or xiii\or xiv\fi\fi}
\def\subsubsubclaimnumber{\ifquiet\else\ifcase\subsubsubclaimno\or a\or b\or 
   c\or d\or e\or f\or g\or \or h\or i\or j\or k\or l\or m\or n\or o\fi\fi}
\def\claimtag{\ifquiet\else
  \ifNoSections
    \ifcase\proofdepth\the\claimno%
    \or\the\claimno.\subclaimnumber
    \or\the\claimno.\subclaimnumber.\subsubclaimnumber
    \or\the\claimno.\subclaimnumber.\subsubclaimnumber
                                                .\subsubsubclaimnumber\fi
  \else
    \ifcase\proofdepth\the\secno.\the\claimno
    \or\the\secno.\the\claimno.\subclaimnumber
    \or\the\secno.\the\claimno.\subclaimnumber.\subsubclaimnumber
    \or\the\secno.\the\claimno.\subclaimnumber.\subsubclaimnumber
                                                .\subsubsubclaimnumber\fi\fi\fi}
\secno=0\claimno=0\proofdepth=0\subclaimno=0\subsubclaimno=0\subsubsubclaimno=0
\NoSectionstrue
\newbox\qedbox
\def\claimname{\ifcase\claimtype Theorem\or Lemma\or Claim\or Corollary\or
               Question\or Definition\or Remark\or Conjecture\fi}
\def\preclaimskip{\removelastskip
    \ifcase\claimtype\goodbreak\vskip 8pt plus 4pt minus 2pt
                  \or\goodbreak\vskip 6pt plus 4pt minus 1pt
                  \or\goodbreak\vskip 5pt plus 4pt minus 1pt
                  \or\goodbreak\vskip 8pt plus 4pt minus 2pt
                  \or\vskip 7pt plus 4pt minus 2pt
                  \or\vskip 7pt plus 4pt minus 2pt
                  \or\vskip 7pt plus 4pt minus 2pt
                  \or\goodbreak\vskip 8pt plus 4pt minus 2pt\fi}
\def\postclaimskip{\ifcase\claimtype         \vskip 4pt plus 2pt minus 2pt
                                          \or\vskip 3pt plus 2pt minus 2pt
                                          \or\vskip 2pt plus 2pt minus 1pt
                                          \or\vskip 4pt plus 2pt minus 2pt
                                          \or\vskip 1pt plus 2pt 
                                          \or\vskip 4pt plus 4pt 
                                          \or\vskip 3pt plus 2pt
                                          \or\vskip 4pt plus 2pt minus 2pt\fi}
\def\claimfont{\ifcase\claimtype
                  \sl\or\sl\or\sl\or\sl\or\sl\or\rm\or\rm\or\sl\fi}
\def\advancetag{\ifcase\proofdepth\advance\claimno by1
                               \or\advance\subclaimno by1
                               \or\advance\subsubclaimno by1
                               \or\advance\subsubsubclaimno by1\fi}
\def\sayclaim#1.#2 #3\par{\ifquiet\else\advancetag\fi
    \preclaimskip\setbox1=\hbox{#1}\setbox2=\hbox{#2}%
    \toks0={#1 }
    \immediate\write16{\ifdim\wd1>0pt\the\toks0
                       \else\claimname\space\fi \claimtag.}%
    \vbox{\noindent
    {\bf\ifdim\wd1=0pt \claimname\else #1\fi\ifquiet.\else\ \claimtag{\ifNoSections.\fi}\fi}%
    \enspace{\ifdim\wd2>0pt\sc #2\enspace\fi}%
    {\claimfont #3\par}}\postclaimskip\quietoff}
\def\theorem{\claimtype=0\sayclaim}

\def\corollary{\claimtype=3\sayclaim}
\def\question{\claimtype=4\sayclaim}

\def\point#1. #2\par{\item{\rm #1.}#2\par}
\def\points#1\cr\par{\medskip\vbox{\let\cr=\point\point#1\par}\par}
\def\df{\it}
\def\prooffont{}
\def\proofsize{}
\def\proofindent{}
\def\proofskip{\badbreak\ifcase\claimtype    \vskip 3pt plus 2pt minus 2pt
                                          \or\vskip 2pt plus 2pt minus 2pt
                                          \or\vskip 1pt plus 2pt minus 1pt
                                          \or\vskip 3pt plus 2pt minus 2pt
                                          \or\vskip 1pt plus 2pt 
                                          \or\vskip 2pt plus 4pt 
                                          \or\vskip 1pt plus 2pt
                                          \or\vskip 3pt plus 2pt minus 2pt\fi}

\def\Goodbreak{\vskip0pt plus.5in\penalty-1000\vskip0pt plus-.5in}
\def\goodbreak{\penalty-500}
\def\badbreak{\penalty500}
\def\Badbreak{\penalty1000}
\def\proof{\message{proof}\removelastskip\Badbreak\proofskip\begingroup
  \advance\proofdepth by1
  \setbox\qedbox=\hbox{\halmos\raise2pt\hbox{\fiverm\claimname}}%
  \prooffont\proofsize\proofindent\noindent{\bf Proof: }}
\def\proofof#1:{\message{proof}\removelastskip\Badbreak\proofskip\begingroup
  \advance\proofdepth by1
  \setbox\qedbox=\hbox{\halmos\raise2pt\hbox{\fiverm#1}}%
  \prooffont\proofsize\proofindent\noindent{\bf Proof of #1: }}
\def\cite[#1]{[{\tenrm{#1}}]\message{[#1]}}
\edef\ref#1{\expandafter\global\expandafter\edef#1{\noexpand\claimtag}}
\newwrite\notes
\openout\notes=\jobname.notes
\long\def\unexpandedwrite#1#2{\def\finwrite{\write#1}%
   {\aftergroup\finwrite\aftergroup{\sanitize#2\endsanity}}}
\def\sanitize{\futurelet\next\sanswitch}
\let\stoken=\space
\def\sanswitch{\ifx\next\endsanity
  \else\ifcat\noexpand\next\stoken\aftergroup\space\let\next=\eat
   \else\ifcat\noexpand\next\bgroup\aftergroup{\let\next=\eat
    \else\ifcat\noexpand\next\egroup\aftergroup}\let\next=\eat
     \else\let\next=\copytoken\fi\fi\fi\fi \next}
\def\eat{\afterassignment\sanitize \let\next= }
\long\def\copytoken#1{\ifcat\noexpand#1\relax\aftergroup\noexpand
  \else\ifcat\noexpand#1\noexpand~\aftergroup\noexpand\fi\fi
  \aftergroup#1\sanitize}
\def\endsanity\endsanity{}

\def\note#1#2{\hbox to2in{\strut#1\quad\dotfill\quad#2}}
\def\boxit#1{\setbox4=\hbox{\kern1pt#1\kern1pt}
  \hbox{\vrule\vbox{\hrule\kern1pt\box4\kern1pt\hrule}\vrule}}
\def\halmos{\hbox{\am\char'3}} 
\def\qed#1\par{\message{.                                }\setbox1=\hbox{#1}%
  \ifdim\wd1>0pt\setbox\qedbox=\hbox{\halmos\raise2pt\hbox{\fiverm #1}}\fi
  \kern5pt\lower 2pt\hbox{\box\qedbox}\proofskip\goodbreak\endgroup}

\def\sectionsymbol{\S}
\def\k{\kappa}
\def\g{\gamma}
\def\a{\alpha}
\def\b{\beta}
\def\d{\delta}
\def\s{\sigma}
\def\t{\tau}
\def\l{\lambda}

\def\I1{\mathop{\hbox{\sc i}_1}}
\def\w{\omega}
\def\P{{\mathchoice{\hbox{\bm P}}{\hbox{\bm P}}
         {\hbox{\tenbm P}}{\hbox{\sevenbm P}}}}

\def\card#1{\left|#1\right|}

\def\Ult{\mathop{\rm Ult}}

\def\elesub{\prec}

\def\unifto{\buildrel\lower 7pt\hbox{$\to$}\over\to}

\def\cof{\mathop{\rm cof}\nolimits}

\def\ran{\mathop{\rm ran}\nolimits}

\def\plus{^{\scriptscriptstyle +}}

\def\in{\mathrel{\mathchoice{\raise 
1pt\hbox{$\scriptstyle\cal\char'62$}}
         {\raise 1pt\hbox{$\scriptstyle\cal\char'62$}}
         {\raise .5pt\hbox{$\scriptscriptstyle\cal\char'62$}}
         {\hbox{$\scriptscriptstyle\cal\char'62$}}}\penalty700{}}
\def\ni{\mathrel{\mathchoice{\raise 1pt\hbox{$\scriptstyle\cal\char'63$}}
                   {\raise 1pt\hbox{$\scriptstyle\cal\char'63$}}
                   {\raise .5pt\hbox{$\scriptscriptstyle\cal\char'63$}}
                   {\hbox{$\scriptscriptstyle\cal\char'63$}}}\penalty700}
\def\of{\mathrel{\mathchoice{\raise 1pt\hbox{$\scriptstyle\subseteq$}}
                   {\raise 1pt\hbox{$\scriptstyle\subseteq$}}
                   {\raise .5pt\hbox{$\scriptscriptstyle\subseteq$}}
                   {\hbox{$\scriptscriptstyle\subseteq$}}}}
\def\fo{\mathrel{\mathchoice{\raise 1pt\hbox{$\scriptstyle\supseteq$}}
                   {\raise 1pt\hbox{$\scriptstyle\supseteq$}}
                   {\raise .5pt\hbox{$\scriptscriptstyle\supseteq$}}
                   {\hbox{$\scriptscriptstyle\supseteq$}}}}
\def\notin{\mathrel{\mathchoice
  {\raise 1pt\hbox{\rlap{$\scriptstyle\;|$}$\scriptstyle\cal\char'62$}}
  {\raise 1pt\hbox{\rlap{$\scriptstyle\kern2pt 
          |$}$\scriptstyle\cal\char'62$}}
  {\raise .5pt\hbox{\rlap{$\scriptscriptstyle\, |$}$\scriptscriptstyle
      \cal\char'62$}}
  {\hbox{\rlap{$\scriptscriptstyle\, |$}$\scriptscriptstyle
     \cal\char'62$}}}%
  \penalty700}

\def\and{\mathrel{\kern1pt\&\kern1pt}}

\def\union{\cup}

\def\compose{\circ}
\def\intersect{\cap}

\def\muchgt{\gg}
\def\lt{\mathrel{\mathchoice{\scriptstyle<}{\scriptstyle<}
   {\scriptscriptstyle<}{\scriptscriptstyle<}}}

\def\[#1]{\left[\vphantom{\bigm|}#1\right]}
\def\<#1>{\langle\,#1\,\rangle}

\def\image{\mathbin{\hbox{\tt\char'42}}}
\def\restrict{\mathbin{\mathchoice{\hbox{\am\char'26}}{\hbox{\am\char'26}}{\hbox{\eightam\char'26}}{\hbox{\sixam\char'26}}}}

\def\beth{\mathord{\hbox{\bm\char'151}}}

\def\st{\mid}
\def\seq<#1>{{\def\st{\mid\penalty650}\left<\,#1\,\right>}}

\def\set#1{\{\,#1\,\}}

\def\th{{\hbox{\fiverm th}}}
\def\ltw{{{\lt}\w}}

\def\lttheta{{\raise 1pt\hbox{$\scriptstyle<$}\theta}}

\def\I1{\mathop{\hbox{\sc i}_1}}

\def\ltl{{{\scriptstyle<}\l}}

\def\hdot{\dot h}

\def\Edot{\dot E}

\font\arrow=line10 scaled \magstep1
\def\makeline#1.{\hbox{\arrow\char#1}}
\def\makearrow#1.#2.{\hbox{\arrow\char#1\llap{\char#2}}}
\def\definelinesandarrows#1.#2.#3.#4.#5.{
   \expandafter\edef\csname#4line\endcsname{\makeline#1.}
   \expandafter\edef\csname#4arrow\endcsname{\makearrow#1.#2.}
   \expandafter\edef\csname#5line\endcsname{\makeline#1.}
   \expandafter\edef\csname#5arrow\endcsname{\makearrow#1.#3.}}
\definelinesandarrows 0.18.9.ne.sw.
\definelinesandarrows 1.21.11.nnne.sssw.
\definelinesandarrows 2.14.13.nnnne.ssssw.
\definelinesandarrows 3.23.15.nnnnne.sssssw.
\definelinesandarrows 4.23.15.nnnnnne.ssssssw.
\definelinesandarrows 10.30.29.nne.ssw.
\definelinesandarrows 16.49.41.neeeeee.swwwwww.
\definelinesandarrows 17.51.43.neeee.swwww.
\definelinesandarrows 19.55.47.nehuh.swhuh.
\definelinesandarrows 24.58.41.neeeeeee.swwwwwww.
\definelinesandarrows 26.62.9.neee.swww.
\definelinesandarrows 33.49.25.neeeee.swwwww.
\definelinesandarrows 35.62.61.nee.sww.
\definelinesandarrows 64.82.73.se.nw.
\definelinesandarrows 65.85.75.ssse.nnnw.
\definelinesandarrows 66.78.77.sssse.nnnnw.
\definelinesandarrows 67.87.79.ssssse.nnnnnw.
\definelinesandarrows 68.87.79.sssssse.nnnnnnw.
\definelinesandarrows 74.94.93.sse.nnw.
\definelinesandarrows 80.113.105.seeeeee.nwwwwww.
\definelinesandarrows 81.115.107.seeee.nwwww.
\definelinesandarrows 99.126.125.see.nww.
\def\sejoin#1#2{\setbox1=\hbox{#1}\setbox2=\hbox{#2}%
  \hbox{\vbox{\hbox{\copy1\kern\wd2}\nointerlineskip
              \hbox{\kern\wd1\box2}}}}
\def\nejoin#1#2{\setbox1=\hbox{#1}\setbox2=\hbox{#2}%
  \hbox{\vbox{\hbox{\kern\wd1\copy2}\nointerlineskip\hbox{\copy1\kern\wd2}}}}
\newdimen\hnudge
\newdimen\vnudge
\newdimen\hnudgedefault
\newdimen\vnudgedefault

\def\SEdefaultnudge{\hnudge=-16pt\vnudge=20pt}
\def\Edefaultnudge{\hnudge=-25pt\vnudge=6pt}
\def\Sdefaultnudge{\hnudge=-8pt\vnudge=20pt}
\def\longEdefaultnudge{\hnudge=-5pt\vnudge=6pt}
\def\nudgeright#1pt{\advance\hnudge by#1pt}
\def\nudgeleft#1pt{\advance\hnudge by-#1pt}
\def\nudgeup#1pt{\advance\vnudge by#1pt}
\def\nudgedown#1pt{\advance\vnudge by-#1pt}
\def\label#1{\smash{\llap{\kern\hnudge
                   \raise\vnudge\rlap{$\scriptstyle#1$}\hfill}}}

\def\SEarrow{\SEdefaultnudge
             \sejoin\seeline{\sejoin\seeline{\sejoin\seeline\seearrow}}}

\def\Sarrow{\Sdefaultnudge\setbox1=\hbox{\SEarrow}
           \hbox{\hskip 10pt\vrule height\ht1\hbox{\arrow\char'77}}}
\def\Earrow{\Edefaultnudge\setbox1=\hbox{\SEarrow}
 \hbox{\raise 2pt\hbox{\vrule height-.4pt depth.8ptwidth\wd1\kern2pt
       \llap{\arrow\char'55}}}}
\def\longEarrow{\longEdefaultnudge\setbox1=\hbox{\SEarrow}
      \rlap{\hskip-1.25\wd1\raise 2pt
            \hbox{\vrule height-.4pt depth.8ptwidth2.5\wd1\kern2pt
            \llap{\arrow\char'55}}}}
\def\trianglediagram#1#2#3#4#5#6{%
    {\def\normalbaselines{\baselineskip0pt\lineskip8pt\lineskiplimit0pt}%
       \matrix{#1& &\cr
               \Sarrow\label{#2}&\SEarrow\label{#3}&\cr
               #4&\Earrow\label{#5}&#6\cr}}}

\QUIET
\def\Edot{\dot E}
\center [submitted to the Proceedings of the American Mathematical Society]

\title Small forcing creates neither\cr
       strong nor Woodin cardinals\cr

\author Joel David Hamkins\cr
	City University of New York\cr
	{\ninett http://www.library.csi.cuny.edu/users/hamkins}\cr

\author  W. Hugh Woodin\cr
            University of California at Berkeley\cr
	{\ninett woodin@math.berkeley.edu}\cr

\abstract. After small forcing, almost every strongness embedding is the lift of a strongness embedding in the ground model. Consequently, small forcing creates neither strong nor Woodin cardinals.

The widely known Levy-Solovay Theorem \cite[LevSol67] asserts that small forcing does not affect the measurability of any cardinal. If a forcing notion $\P$ has size less than $\k$, then $\k$ is measurable in $V^\P$ if and only if it is measurable in $V$; every measure on $\k$ in $V$ extends uniquely to a measure on $\k$ in $V^\P$; the corresponding embeddings lift uniquely from $V$ to $V^\P$, and all the measures in $V^\P$ arise in this way. The argument generalizes to show the same fact for small forcing with the other large cardinals which are witnessed by the existence of certain kinds of ultrapowers, such as strongly compact cardinals, supercompact cardinals, almost huge cardinals, huge cardinals, and so on. Missing, however, from the scope of the theorem are both the strong cardinals and the Woodin cardinals, whose embeddings are not simple ultrapowers, but directed systems of them. Several years ago the second author of this paper successfully treated these extender embeddings and proved that small forcing creates neither strong nor Woodin cardinals. Since his proof seems not to be widely known, we would like to present it here along with a level-by-level version which aims to obtain the result for cardinals which are only partially strong. 

\theorem Main Theorem. After small forcing, a cardinal $\k$ is strong if and only if it was strong in the ground model.

\theorem Level-by-level Version. If $g\of\P$ is $V$-generic for the small forcing $\P$, then for every ordinal $\l$ (except possibly the limit ordinals with $\cof(\l)\leq\card{\P}\plus$), every natural $\l$-strongness embedding $j:V[g]\to M[g]$ in the extension lifts a $\l$-strongness embedding $j\restrict V:V\to M$ definable in the ground model. 

\corollary. After small forcing, a cardinal $\theta$ is Woodin if and only if it was Woodin in the ground model.

Let us define the relevant concepts. A forcing notion $\P$ is {\df small} relative to $\k$ when $\card{\P}<\k$. A cardinal $\k$ is {\df strong} when for every $\l$ it is {\df $\l$-strong}, so that there is a $\l$-strongness embedding $j:V\to M$, one with critical point $\k$ and $V_\l\of M$. A cardinal $\theta$ is a {\df Woodin} cardinal when for every $A\of\theta$ there is a $\k<\theta$ which is $\lttheta$-strong for $A$, meaning that for every $\l<\theta$ there is an embedding $j:V\to M$ with critical point $\k$ such that $j(A)\intersect\l=A\intersect\l$; such embeddings are said to be {\df $\l$-strong for $A$}. By coding information into the set $A$, one can always find such a $j$ which is $\l$-strong. A $\l$-strongness embedding $j:V\to M$ is {\df natural} when $M=\set{j(h)(s)\st h\in V\and s\in\g^\ltw}$, where $\g=\card{V_\l}$. Equivalently, we could require that $M=\set{j(h)(s)\st h\in V\and s\in V_\l}$. 

Before taking up the theorem, let us prove a few general facts about strongness embeddings. First, we claim that every strongness embedding factors through a natural embedding. To see this for any particular $\l$-strongness embedding $j:V\to M$, where $\l>\k$, let $X=\set{j(h)(s)\st h\in V\and s\in\g^\ltw}$. By checking the Tarski-Vaught criterion, it is easy to see that $X$ is an elementary substructure of $M$ covering $\ran(j)$, and so if $\pi:X\to M_0$ is the Mostowski collapse of $X$, we obtain a map $j_0:V\to M_0$ by simply defining $j_0=\pi\compose j$. Thus also $j=k\compose j_0$, where $k=\pi^{-1}$, and so $j$ factors through $j_0$. $$\trianglediagram{V}{j_0}{j}{M_0}{k}{M}$$ Since $V_\l\of X$ it follows that $\pi\image V_\l=V_\l$ and $V_\l\of M_0$, and since $M_0=\ran(\pi)$ it follows that $M_0=\set{j_0(h)(s)\st h\in V\and s\in \g^\ltw}$. Thus, $j_0$ is a natural $\l$-strongness embedding, as desired. Note that if $j$ is $\l$-strong for $A$, then so also is $j_0$, because the critical point of $k$ is at least $\l$. The point is therefore that when $\theta$ is a Woodin cardinal, one can take the witnessing embeddings to be natural. 

Second, while strongness embeddings in general need not satisfy any closure properties, the natural strongness embeddings are much better behaved. Specifically, we claim that if $j:V\to M$ is a natural $\l$-strongness embedding with critical point $\k$ and $\l$ is either a successor ordinal or a limit ordinal of cofinality above $\k$, then $M$ is closed under $\k$-sequences. Otherwise, $M$ is closed under ${\scriptstyle <}\cof(\l)$-sequences. To see this, we may suppose that $\l>\k$ and $M=\set{j(h)(s)\st h\in V\and s\in V_\l}$. In the first case, suppose that $\l=\xi+1$ and $\<j(h_\a)(s_\a)\st\a<\k>$ is a $\k$-sequence of elements from $M$, with each $s_\a\in V_{\xi+1}$. Since a $\k$-sequence of subsets of $V_\xi$ can be coded with a single subset of $V_\xi$, it follows that $\<s_\a\st\a<\k>$ is in $M$. Then, since the sequence $\<j(h_\a)\st\a<\k>=j(\<h_\a\st\a<\k>)\restrict\k$ is in $M$, it follows that $\<j(h_\a)(s_\a)\st\a<\k>$ is in $M$, as desired. For the next case, when $\l$ is a limit ordinal of cofinality larger than $\k$, then on cofinality grounds the sequence $\<s_\a\st\a<\k>$ is in $V_\l$, and hence in $M$, so again $\<j(h_\a)(s_\a)\st\a<\k>$ is in $M$, as desired. Finally, if $\l$ is a limit ordinal, $\b<\cof(\l)\leq\k$ and $\<j(h_\a)(s_\a)\st\a<\b>$ is a sequence of elements from $M$, then again on cofinality grounds we know $\<s_\a\st\a<\b>$ is in $V_\l$ and hence in $M$, and so $\<j(h_\a)(s_\a)\st\a<\b>$ is in $M$, as desired. 

Now we are ready to prove the main theorem. 

\theorem Main Theorem. Suppose $g\of\P$ is $V$-generic for forcing $\P$ of size less than $\k$. Then $\k$ is strong in $V$ if and only if it is strong in $V[g]$.

\proof We may assume that $\P\in V_\k$. Let $\d=\card{\P}<\k$. The forward direction of the theorem is trivial because any $\l$-strongness embedding $j:V\to M$ in $V$ lifts to an embedding $j:V[g]\to M[g]$ in $V[g]$ by simply defining $j(\t_g)=j(\t)_g$ for any name $\t\in V$. This embedding witnesses the $\l$-strongness of $\k$ in $V[g]$ because $V[g]_\l=V_\l[g]$. For the converse, suppose that $\k$ is $\l$-strong in $V[g]$ with natural embedding $j:V[g]\to M[g]$, where $\l>\k$ is either a successor ordinal or a limit ordinal with $\cof(\l)>\d\plus$. To prove the theorem, we will show that $j\restrict V:V\to M$ is definable in $V$ and witnesses there that $\k$ is $\l$-strong. Since $j$ is natural, we know that $M[g]=\set{j(h)(s)\st h\in V[g]\and s\in \g^\ltw}$, where $\g=\beth_\l$. It is easy to verify that $V_\k=M_\k$ because sets with rank below the critical point are fixed by $j$. Moreover, if $A\of\k$ in $V$ then $A=j(A)\intersect\k$ and so $P(\k)^V\of M$. Conversely if $A\of\k$ in $M$ then every initial segment of $A$ is in $V$ and so $A$ must also be in $V$ since the forcing is $\k$-c.c. It follows that $V_{\k+1}=M_{\k+1}$. Let $E$ be the induced $V$-extender, namely, $E=\set{\<A,s>\st s\in\g^\ltw\and A\in V_{\k+1}\and s\in j(A)}$. 

We claim that the restricted embedding $j\restrict V:V\to M$ is precisely the embedding induced by the extender $E$. It suffices to show that $M=\set{j(h)(s)\st h\in V\and s\in\g^\ltw}$, for then the map which takes the $E$-equivalence class of $\<h,s>$ in $\Ult(V,E)$ to $j(h)(s)$ will be an isomorphism. By the assumption on $j$ we know that any set $a$ in $M$ has the form $j(h)(s)$ for some function $h:\k^n\to V$ in $V[g]$ and some $s\in\g^n$ and $n\in\w$; the difficulty is to find such a function $h$ in $V$. But since $h=\hdot_g$ for some name $\hdot\in V$ and $j(\hdot)_g(s)=a$ there must be a condition $p\in g$ which forces over $M$ that $j(\hdot)(\check s)=\check a$. In $V$ define the function $\bar h(x)=y$ when $p$ forces that $\hdot(\check x)=\check y$. It follows that $j(\bar h)(s)=a$, and so $a$ has the desired form. 

The point now is that in order to show that $j\restrict V$ is definable (from parameters) in $V$, it suffices to show that $E$ is in $V$. For this, we will first prove the claim that if $F\of E$ and $\card{F}=\d$ then there is $F^*\in V$ such that $F\of F^*\of E$. Fix such an $F$, and let $\s$ be the set of generators which appear in $F$, that is, the set of elements appearing in the various $s$ which appear in $F$. So $\s\of \g$ and $\card{\s}\leq\d$. Consequently, by the hypothesis on $\l$ and the remarks on the closure of natural strongness embeddings, $\s$ is in $M[g]$. So the set $\s$ has names in both $M$ and $V$. By iteratively using such names from $V$ and $M$, we may construct an increasing $\d\plus$ sequence of sets $\vec\s=\<\s_\a\st\a<\d\plus>$, beginning with $\s_0=\s$, such that every $\s_\a$ has cardinality $\d$, for cofinally many $\a$ the set $\s_\a$ is in $V$ and for cofinally many $\a$ it is in $M$. Let $\t=\union_\a\s_\a$. Since $\vec\s\in M[g]$ by the closure of the embedding, the sequence has names in both $M$ and $V$. In $V$ there are conditions in $\P$ which force the various $\s_\a$ to appear and so since $\d\plus$ is regular there must be a single condition which decides unboundedly many of the $\dot\s_\a$ in $\dot{\vec\s}$. Thus, $\t\in V$. Similarly, using a name in $M$ for $\vec\s$ there must be a single condition deciding unboundedly many of the $\s_\a$, and so $\t\in M$. Let $\mu=\set{A\in V_{\k+1}\st \t\in j(A)}$. Since this is a $V$-measure in $V[g]$, the Levy-Solovay Theorem for measurable cardinals tells us that it lies in $V$. From $\mu$ we can compute the set $F^*=\set{\<A,s>\st s\in\t^\ltw\and A\in V_{\k+1}\and s\in j(A)}=E\restrict\t$, since any such $s$ corresponds to a certain projection of $\t$. That is, if $s$ consists of the $\a_0^\th, \ldots, \a_n^\th$ elements of $\t$, then $\<A,s>\in F^*$ exactly when $\pi^{-1}A\in\mu$ where $\pi(t)$ restricts $t$ to its $\a_0^\th, \ldots, \a_n^\th$ elements. Consequently, $F^*$ is in $V$, as we claimed. 

Now let us show that $E\in V$. Since $E\in V[g]$ it has a name $\Edot$ in $V$. Suppose $\eta\muchgt\gamma$ and let $X\elesub V_\eta$ be an elementary substructure of size $\d$ which contains $\P$, $F^*$, $\Edot$ and every element of $\P$. It follows that $g$ is $X$-generic, that $X[g]\elesub V_\l[g]$ and that $E\in X[g]$. Let $F=E\intersect X[g]$. By the previous paragraph, there is a set $F^*\in V$ such that $F\of F^*\of E$. Since $F=F\intersect X\of F^*\intersect X\of E\intersect X=F$, it follows that $F=F^*\intersect X$ is in $V$. Thus, there must be a condition $p\in g$ which decides ``$\<A,s>\in\Edot\,$'' for all $\<A,s>$ is in $X$. By elementarity, therefore, $p$ decides ``$\<A,s>\in\Edot\,$'' for all $\<A,s>$ in $V$. So $E$ lies in $V$. 

We have argued that the restricted embedding $j:V\to M$ is therefore definable from $E$ in $V$. In order to finish the proof of the theorem, then, it suffices for us to show that $V_\l\of M$. Certainly $V_\l\of M[g]$ since the full embedding was $\l$-strong, and so this must be forced by some condition $p\in g$. That is, $V_\l\of M[g]$ for any generic below $p$. Thus, for any $a\in V_\l$ there is a name $\dot a\in M$ such that $p$ forces (over $V$) that $\check a=\dot a$. It follows using $p$ in $M$ that $a\in M$. That is, $V_\l\of M$. Thus, the restricted embedding $j:V\to M$ witnesses that $\k$ is $\l$-strong in $V$.\qed

The proof actually established the level-by-level version that we stated in the beginning of this paper: 

\theorem Level-by-level Version. If $g\of\P$ is $V$-generic for the small forcing $\P$, then for every ordinal $\l$ (except possibly the limit ordinals with $\cof(\l)\leq\card{\P}\plus$), every natural $\l$-strongness embedding $j:V[g]\to M[g]$ in the extension lifts a $\l$-strongness embedding $j\restrict V:V\to M$ definable in the ground model. 

We are left with the intriguing possibility that small forcing could make a $\ltl$-strong cardinal $\l$-strong when $\l$ is a limit ordinal with small cofinality. Certainly it seems possible that if there are several kinds of $\l$-extenders in the ground model, then after adding a Cohen real, say, one could consult the new real $x$ and glue together various pieces of these extenders in such a way that the resulting embedding $j:V[x]\to M[x]$ could not be the lift of an embedding from the ground model; the idea being that from the restricted embedding $j:V\to M$ one would be able to recover the real $x$. In such a situation, the strict analogue of the Levy-Solovay theorem for strongness embeddings would fail, because not every natural strongness embedding  would lift an embedding from the ground model. More difficult, however, is the question of whether after small forcing one could glue together various $\ltl$-strongness extenders to make a $\l$-strong embedding in the extension even when no $\l$-strong embedding exists in the ground model. That is, the following question remains open:

\question. Can small forcing increase the degree of strongness of a cardinal?

Of course, the only possibility for an affirmative answer is that a $\ltl$-strong cardinal $\k$ is made $\l$-strong in the extension for some limit ordinal of small cofinality. 

Let us now treat the case of Woodin cardinals. 

\corollary. After small forcing, a cardinal is Woodin if and only if it was Woodin in the ground model. 

\proof As in the theorem, a trivial lifting argument shows that small forcing cannot destroy any Woodin cardinal. Conversely, suppose that $\theta$ is Woodin in $V[g]$, a small forcing extension, and that $A\of\theta$ in $V$. There must be a cardinal $\k<\theta$ such that for every $\l<\theta$ there is a natural $\l$-strongness embedding $j:V[g]\to M[g]$ with critical point $\k$ such that $j(A)\intersect\l=A\intersect\l$. If $\l$ is a successor ordinal, then by the Level-by-Level version of the Main Theorem the restricted embedding $j:V\to M$ satisfies the same hypothesis in $V$. Thus, since the successor ordinals are unbounded in $\theta$, this means that $\theta$ is Woodin in $V$, as desired.\qed

\quiet\section Bibliography

\nopagenumbers
\parindent=0pt
\newbox\Article
\newbox\Journal
\newbox\Author
\newbox\Vol
\newbox\No
\newbox\Year
\newbox\Page
\newbox\Book
\newbox\Publisher
\newbox\Pubaddr
\newbox\Key
\newbox\Editor
\newbox\Comment
\newbox\Note
\def\entry#1#2\par{\item{#1\quad}\hskip-1.1em#2\par}
\def\article#1{\setbox\Article=\hbox{\sl #1, }}
\def\journal#1{\setbox\Journal=\hbox{\rm #1 }}
\def\author#1{\setbox\Author=\hbox{\sc #1, }}
\def\vol#1{\setbox\Vol=\hbox{\bf #1 }}
\def\no#1{\setbox\No=\hbox{no. #1 }}
\def\year#1{\setbox\Year=\hbox{\rm({\oldstyle #1}) }}
\def\page#1{\setbox\Page=\hbox{\rm p. #1 }}
\def\book#1{\setbox\Book=\hbox{\it #1, }}
\def\publisher#1{\setbox\Publisher=\hbox{\rm #1, }}
\def\pubaddr#1{\setbox\Pubaddr=\hbox{\rm #1, }}
\def\key#1{\setbox\Key=\hbox{#1}}
\def\editor#1{\setbox\Editor=\hbox{\rm(#1, Ed.) }}
\def\comment#1{\setbox\Comment=\hbox{\rm #1}}
\def\note#1{\setbox\Note=\hbox{\rm #1 }}
\def\ref#1\par{\smallskip{#1
  \entry{\ifhbox\Key\unhbox\Key\else[\ ]\fi}%
  \unhbox\Author\unhbox\Note
  \ifhbox\Book \unhbox\Book\unhbox\Publisher\unhbox\Pubaddr
               \unhbox\Editor\unhbox\Page\unhbox\Year\unhbox\Comment
  \else \unhbox\Article\unhbox\Journal\unhbox\Vol\unhbox\No\unhbox\Editor
        \unhbox\Page\unhbox\Year\unhbox\Comment\fi\par}}

\tenpoint\tightlineskip

\ref
\author{Azriel Levy, Robert M. Solovay}
\article{Measurable cardinals and the Continuum Hypothesis}
\journal{IJM}
\vol{5}
\year{1967}
\page{234-248}
\key{[LevSol67]}

\bye